\documentclass[conference]{IEEEtran}
\IEEEoverridecommandlockouts

\usepackage{amsmath,amssymb,amsfonts}

\usepackage{graphicx}
\usepackage{textcomp}
\usepackage{xcolor}
\usepackage{algorithm}
\usepackage{algpseudocode}
\usepackage[
    backend=biber,
    citestyle=numeric,
    bibstyle=numeric,
    maxcitenames=7,
    style=numeric,
    maxbibnames=99,
    sorting=none
]{biblatex}
\DeclareNameAlias{sortname}{family-given}
\usepackage{tikz}
\def\BibTeX{{\rm B\kern-.05em{\sc i\kern-.025em b}\kern-.08em
    T\kern-.1667em\lower.7ex\hbox{E}\kern-.125emX}}
\begin{document}

\title{Space Mapping Optimization using Neural Networks for Efficient Parameter Estimation\\}

\author{\IEEEauthorblockN{Dhruvil Kamleshkumar Kotecha}
\IEEEauthorblockA{
Hamilton, Canada \\}
}

\maketitle

\begin{abstract}
This project focuses on optimizing input parameters of a partial derivative function of a fine model using Neural network-based Space Mapping Optimization (SMO). The fine model is known for its high accuracy but is computationally expensive. On the other hand, the coarse model is represented by a neural network, which is much faster but less accurate. The SMO approach is applied to bridge the gap between these two models and estimate the optimal input parameters for the fine model. Additionally, this project involves a comprehensive review of previously available Neuro Modeling Space Mapping techniques, which are also used in this project to enhance the optimization process. By utilizing SMO with a neural network-based coarse model, we aim to demonstrate the effectiveness of this method in optimizing complex functions efficiently. The proposed approach of using Neural Network based Space Mapping offers a promising solution to this optimization problem..
\end{abstract}

\section{\textbf{Introduction}}
\vspace{3pt}
Optimizing complex models has become increasingly crucial in various fields, such as aerospace, automotive, and civil engineering. Optimization techniques are utilized to enhance the design of products and processes, leading to more efficient and cost-effective solutions. Among these techniques, Space Mapping Optimization (SMO) has shown great potential in optimizing the input parameters of a partial derivative function of a fine model. However, optimizing fine models can be computationally expensive, presenting a significant challenge in the optimization process. In such cases, a coarse model represented by a neural network is often used, which is faster but less accurate.

This project aims to explore the use of neural network-based Space Mapping Optimization (SMO) to enhance the optimization of a partial derivative function of a fine model. The objective is to efficiently estimate the optimal input parameters for the fine model by bridging the gap between a fine model and a neural network-based coarse model. The use of neural networks to construct surrogate models for complex physical systems is an efficient model, which can replace computationally expensive fine models while still providing accurate approximations of the system's behavior. Neurospace modeling can be incorporated into the Space Mapping Optimization (SMO) technique to enhance the optimization process's efficiency. By utilizing neural networks to represent the coarse model, the computational cost of the optimization process can be reduced while maintaining high accuracy. Furthermore, the use of Neurospace modeling can help overcome challenges associated with traditional surrogate modeling approaches, such as the curse of dimensionality and the need for extensive training data.

The project's primary objective is to demonstrate the effectiveness of using SMO with a neural network-based coarse model to optimize complex functions efficiently.
\section{\textbf{Background}}
\vspace{3pt}
\subsection{Space Mapping Optimization}
\vspace{2pt}
Space mapping optimization is a well-established optimization technique that is utilized to optimize computationally expensive and complex simulations. Initially developed for optimizing microwave circuits in microwave engineering, this technique has been effectively applied to diverse fields such as electromagnetics, fluid mechanics, and structural optimization. Numerous researchers have contributed to this area, including Bandler et al. (2001)\cite{bakr} who introduced the space mapping technique.

The fundamental concept of space mapping optimization involves creating a coarse model that is computationally inexpensive and easy to evaluate. The design parameters of the fine model are then mapped onto the coarse model, and the coarse model is optimized utilizing standard optimization methods. The optimization outcomes are then mapped back to the fine model, and the process is iterated until convergence is achieved.

One of the primary benefits of space mapping optimization is its ability to significantly reduce the computational time and cost necessary for design optimization. This is especially useful in cases where the fine model requires hours or days to simulate, making it impractical for optimization. By utilizing a coarse model, the optimization process can be accelerated significantly, particularly when dealing with large-scale optimization problems, such as designing complex systems or structures.

The design parameters of the coarse and fine models are represented by vectors $x_{c}$ and $x_{f}$, respectively, and their corresponding model responses are represented by $R_{c}(x_{c})$ and $R_{f}(x_{f})$. Notably, $R_{c}$ can be computed more rapidly than $R_{f}$, but $R_{f}$ produces more precise results.

The primary objective of space mapping optimization is to establish a suitable mapping function, represented by P, that maps the fine model parameter space, $x_{f}$, to the coarse model parameter space, $x_{c}$. The mapping function can be expressed as follows \cite{bakr2000review}:
\begin{equation}
x_{c}=P(x_{f})\label{eq}
\end{equation}

The objective is to find a valid mapping function P that satisfies the following condition:

\begin{equation}
R_{c}(P(x_{f}))\approx R_{f}(x_{f}).\label{eq1}
\end{equation}

Once a valid mapping function P has been determined for the region of interest, the coarse model can be used to carry out rapid and precise simulations within the same region. $P$ acts as a link between the coarse and fine models, allowing the optimization process to proceed more effectively.

The optimization process starts by creating the coarse model and training it with appropriate data. The fine model is then defined, and $P$ is initialized. The optimization is performed iteratively until convergence is reached. During each iteration, $P$ is utilized to optimize the coarse model, and the optimization outcomes are mapped back to the fine model. This process is continued until the optimized response of the coarse model is sufficiently close to that of the fine model. Once $P$ has been identified, it can be used to optimize the design parameters of the fine model in a more efficient manner.

\vspace{2pt}
\subsection{\textbf{Use of ANN}}
\vspace{3pt}

Artificial neural networks (ANNs) are a powerful technique for approximating measurable functions to a desired level of accuracy. When a deterministic relationship exists between input and output variables, multi-layer feed-forward neural networks can achieve this with superior computational efficiency and accuracy, especially in modeling microwave circuit yield and statistical design.

Studies by Zaabab et al. (1995)\cite{Zaabab1995ANN} and Burrascano et al. (1998)\cite{burrascano1998neural} have demonstrated the benefits of using ANNs to optimize the yield of microwave circuits. ANNs are capable of modeling complex relationships and nonlinearities that other approaches may not capture accurately.

Training an ANN involves modifying the internal parameters of the neural network to create an ANN model that best fits the training data. This process involves employing optimization techniques such as backpropagation to update the weights in the network. By identifying the optimal combination of weights and biases, the ANN can accurately map inputs to outputs, making it a powerful tool for prediction and modeling within a designated region of interest.

It is important to carefully consider the complexity of the ANN, as overly simple networks may result in underfitting, and overly complex networks may result in overfitting. Choosing an appropriate level of complexity is crucial to achieving accurate results.

The ultimate objective of using ANNs is to determine optimal values for the internal parameters of the model, such that the response of the coarse model closely approximates the response of the fine model at all learning points\cite{phdthesis}.

The response of the fine model be $a_{i}=R_{f}(x_{f_{i}})$ and the response from the coarse model be $b_{i}=R_{c}(P(x_{f_{i}}))$. The coarse model( Neural Network in this case) is trained using the same inputs as the fine model. Optimization is performed until the following difference(loss function)equation (3) holds true. $Note:$ the value of the threshold $\epsilon$ should be small enough such that a close approximation between the response of optimized coarse model and fine model is guaranteed.

\begin{equation}
    l(\underline{a},\underline{b}) \leq \epsilon
\end{equation}
where 
\begin{equation}
    \underline{a}=\{R_{f}(x_{f_{1}}),R_{f}(x_{f_{2}}),\dots,R_{f}(x_{f_{k}})\}
\end{equation}
\begin{equation}
    \underline{b}=\{R_{c}(P(x_{f_{1}})),R_{c}(P(x_{f_{2}})),\dots,R_{c}(P(x_{f_{k}}))\}
\end{equation}
\begin{equation}
    l(a_{i},b_{i}) = a_{i}-b_{i}, \forall i \in \{1,2,\dots,k\} 
\end{equation}
loss function $l$ can be a simple subtraction formula. And $k$ is the number of training examples. To find the mapping, the following optimization problem\cite{bakr123} needs to be solved
\begin{equation}
    \min_{w} \|[l(a_{i},b_{i}^{T})]^{T}\|,\quad \forall i \in \{1,2,\dots,k\}
\end{equation}
The vector $w$ represents a set of internal parameters that are essential for the effective functioning of a neural network. These parameters are typically used to define the weights and bias associated with the various connections between the different neurons within the network architecture. 

\subsection{Problem Formulation}
\vspace{3pt}
One of the advantageous uses of space mapping optimization along with neural network is solving the partial differential equations. A fine model solver can be used for this purpose.  Using ANN  as a coarse model ( even if it is poorly trained initially) can be an effective way that will eventually lead to highly accurate and less expensive solution. To give an example, this project involves solving a heat transfer problem by numerically solving a complex differential equation. The problem involves the transfer of heat through a medium with certain properties and boundary conditions. The problem is formally defined as following:

\begin{figure}
    \centering

    \includegraphics[width=0.6\textwidth]{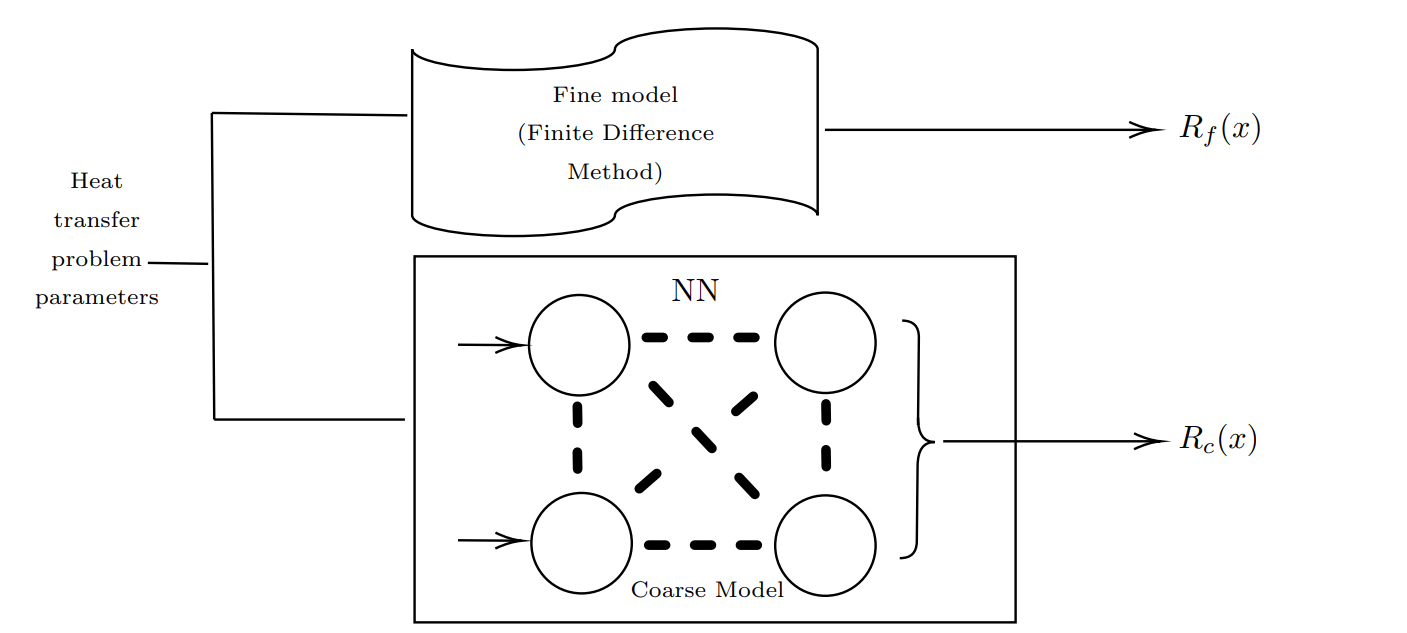}
    \caption{Problem Formulation}
    \label{fig:my_label}
\end{figure}

In this project, a fine model solver, such as a finite element method or a finite difference method, can be used to solve the heat transfer differential equation. The fine model will provide accurate results, but it can be computationally expensive. Therefore, a neural network can be used as a surrogate for the fine model to reduce the computational cost.

The first equation is a differential equation, which describes how a physical quantity (in this case, the temperature) changes over time and space. The equation is:
\begin{equation}
    \frac{\partial T}{\partial t} = \alpha \nabla^2 T + S(x,t)
\end{equation}
where $\frac{\partial T}{\partial t}$ is the rate of change of temperature with respect to time, $\alpha$ is a thermal diffusivity constant, $\nabla^2 T$ is the Laplacian operator applied to the temperature field (which describes how the temperature changes over space), and $S(x,t)$ is a source term that describes any external heat sources or sinks.

The Laplacian operator is a mathematical operator that describes how a function (in this case, the temperature field) changes over space. It is defined as:
\begin{equation}
    \nabla^2 T = \frac{\partial^2 T}{\partial x^2} + \frac{\partial^2 T}{\partial y^2} + \frac{\partial^2 T}{\partial z^2}
\end{equation}
where $\frac{\partial^2 T}{\partial x^2}$, $\frac{\partial^2 T}{\partial y^2}$, and $\frac{\partial^2 T}{\partial z^2}$ are the second partial derivatives of temperature with respect to the x, y, and z directions, respectively.

The second equation is an algebraic equation that relates the heat flux at the boundary of the object to the temperature at the boundary. It is called the Robin boundary condition and is given by:
\begin{equation}
    k \frac{\partial T}{\partial n} = h(T-T_\infty)
\end{equation}
where $k$ is the thermal conductivity of the material, $\frac{\partial T}{\partial n}$ is the normal derivative of temperature at the boundary (which describes how the temperature changes in the direction perpendicular to the boundary), $h$ is the heat transfer coefficient, $T_\infty$ is the temperature of the surrounding fluid (which is assumed to be constant), and $T$ is the temperature at the boundary.

As stated previously, the coarse model is a poorly trained neural network and fine model is the finite difference method. Finite difference method is a numerical method for solving differential equations by approximating the derivatives using finite differences. In this method, the domain is discretized into a grid of points, and the values of the function (in this case, the temperature) at each point are calculated based on the values at neighboring points. 

The finite difference method is widely used for solving partial differential equations like the one described above. The finite difference method involves three methods of approximation: forward difference, backward difference, and central difference. It is used to obtain a numerical solution to a partial differential equation in a bounded domain. In this method, the solution to the PDE is replaced with an approximation using a finite number of points in the domain. The accuracy of the numerical solution generally increases as the number of points increases \cite{book}. 

The explicit method \cite{parkinson2022solving} for solving PDEs known as forward difference utilizes the current point $x_i$ and the next grid point to approximate $f_i$.
\begin{center}
    $f_{i}=f(x_{i+1})-f(x_{i})$
\end{center}
On the other hand, backward process uses the current and the previous point
\begin{center}
    $f_{i}=f(x_{i-1})-f(x_{i})$
\end{center}
And lastly, the central difference uses the previous point point and next point from the current point.
\begin{center}
    $f_{i}=f(x_{i+1})-f(x_{i-1})$
\end{center}
This project uses central difference scheme is used to approximate the second-order partial derivative of T with respect to x. The finite difference scheme is then used to discretize the Laplacian operator in the partial differential equation, which results in a system of algebraic equations that can be solved numerically. This method discretizes the domain and solves the given partial differential equation in space, and then it uses the trapezoidal rule method to integrate the solution in time.

\begin{figure}
    \centering
    \begin{tikzpicture}
\draw (0,0) node[anchor=east] {$i-1,n$} -- (4,0) node[anchor=west] {$i+1,n$};

\coordinate (M) at (2,0);

\draw (M) -- (2,2) node[anchor=south] {$i,n+1$};
\draw (M) -- (2,0) node[anchor=north] {$i,n$};
\end{tikzpicture}
    \caption{The finite difference method for the heat equation}
    \label{fig:my_label}
\end{figure}
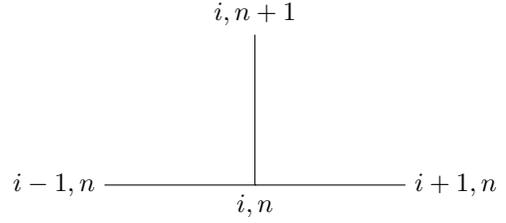

The ultimate goal of this project is to find optimized parameters $\alpha$, $S(x,t)$, $k$, $h$, $T_\infty$ using the space mapping optimization. This project employs a multi-layer neural network, or deep network, which includes at least two hidden layers in addition to an input and output layer. This is in contrast to a shallow network, which only has one hidden layer, one input layer, and one output layer. Shallow networks have a direct mapping from the input to output, which is often insufficient for many use cases and not as efficient as a deep network with multiple hidden layers
. According to Telgarsky's 2016 paper on the ``Benefits of Depth in Neural Networks," there exist neural networks (NNs) with $O(k^{3})$ layers, with a constant number of nodes per layer, and a constant number of distinct parameters. These networks cannot be effectively approximated by networks with O(k) layers, unless the latter are of exponential size. This holds true for networks that feature semi-algebraic gates, including ReLU and piece-wise polynomial functions.
\section{Methodology}
\vspace{3pt}
This project employs a partial differential equation as the underlying model and utilizes a numerical solver as the fine model, while a lightly trained neural network serves as the coarse model. Specifically, the project trains the neural network to predict the solution of the differential equation based on input parameters, thereby emulating the fine model's response, which involves computationally expensive numerical solving techniques. This project follows the following algorithm:
\begin{algorithm}
\caption{}\label{alg:cap}
\begin{algorithmic}[1]
\State Define Partial Derivative function $f_{d}(\underline{x})$ to be solved
\State Build the required neural network (coarse model)
\State Train the neural network with training data
\State Define the fine model $R_{f}$
\State Initialize i=1
\State Initialize threshold $\epsilon$
\While{$(R_{f}(x_{f_{i}})-R_{c}(x_{c^{*}})\leq \epsilon)$}
\State Optimization on coarse model, find $x_{c^{*}}$
\State Refine the solution using fine model $R_{f}(x_{f})$
\State Increase the coarse model with $x_{f}$
\EndWhile
\end{algorithmic}
\end{algorithm}

$Step:1$ defines the problem of heat transfer problem by providing the input parameters to the problem which need to be optimized. The fine model is a complex differential equation that describes the behavior of a heat transfer problem. $Step:2$ requires to build and define the architecture of the neural network, including the number of layers, number of neurons in each layer, and activation functions. This project uses 2 hidden layers with Rectified Linear activation Unit(ReLU) as the activation function for the hidden layers.\begin{center}
     \[   
    ReLU(z) = max(0,z)=
     \begin{cases}
       \text{z,} &\quad\text{if z}> 0,\\
       \text{0,} &\quad\text{otherwise} \\
     \end{cases}
\]

\end{center} And the backpropogation algorithm uses the gradient of the activation in the process of gradient descent while updating the weights of the neural network. The $gradient(ReLU(z))$ is not defined at $z=0$, however, for the convenience I have considered the following.\begin{center}
            \[   
    gradient(ReLU(z)) =
     \begin{cases}
       \text{1,} &\quad\text{if z}\ge 0,\\
       \text{0,} &\quad\text{if z}<0 \\ 
     \end{cases}
\]
\end{center}$Step:3$ trains the neural network.The neural network takes the same set of input parameters as the fine model and predicts the output value (the temperature distribution) based on those inputs. The training data for the neural network is generated by running the fine model for a set of random input parameter values and recording the resulting output values.\break$Step:4$\  defines the fine model function which takes in a set of parameters including the diffusion coefficient alpha, the source term, the thermal conductivity, the heat transfer coefficient, the temperature at infinity, the spatial coordinate, time, and temperature at the current point. The function returns the rate of change of temperature with respect to time based on the differential equation. \break$Step:6$ involves setting an accuracy threshold $\epsilon$, which represents the level of error that is acceptable in the final solution. The while loop will continue until the difference between the predictions of the fine model and the coarse model for a specific input is less than or equal to the specified accuracy threshold. \break$Step:7$ is the main loop of the algorithm. It continues until the difference between the predictions of the fine model and the coarse model for a specific input is less than or equal to the specified accuracy threshold. If this condition is not met, the loop will continue to execute. \break$Step:8$ performs the optimization to find $x_{c}^{*}$. This project uses conjugate gradient method[*] to perform optimization for the coarse model. \break$Step:9$ Once the optimization is performed,the algorithm then uses the optimized values of the parameters($x_{c}^{*})$ to solve the fine model. \break$Step:10$ the fine model solution is then used to update the neural network with new training examples.

To put everything together, the algorithm outlines a process for solving partial derivative functions using a combination of a coarse neural network model and a fine model. The coarse model is trained using training data, and optimization is performed on this model to find an approximate solution. This solution is then refined using the fine model, and the coarse model is updated with this refined solution. The process continues until a desired level of accuracy is achieved. 

\section{\textbf{Analysis}}
\vspace{3pt}
\subsection{Result}
\vspace{3pt}
For the coarse model optimization, there many methods that can be used like Conjugate Gradient (CG) Method \cite{Shewchuk1994AnIT}, Broyden-Fletcher-Goldfarb-Shanno (BFGS) Method and Nelder-Mead Method \cite{Nelder1965ASM}. This project worked on the $\underline{x}_{f}$ values as the training data, where $\underline{x}_{f}$ = \{$\alpha$, $S(x,t)$, $k$, $h$, $T_\infty$\} and the random values of $x,t$ (used in equation 8) were assigned and their appropriate $T$ values were computed using the analytic method as shown in the following equation. It directly calculates the temperature at any given point and time using the analytical solution of the heat equation. This analytical solution method was found in the book by \emph{Poirier and Geiger}\cite{Poirier2016}\begin{equation}
\footnotesize
    T = \sin(\pi x) \exp\left(-\frac{k}{h^2} \pi^2 t\right) + T_{\infty} + \frac{S(x,t)}{k}\left(1 - \exp\left(-\frac{k}{h^2} \pi^2 t\right)\right)
\end{equation}These data was then used to train the artificial neural network (coarse model). Algorithm [1] was then applied and the solution to the heat transfer partial derivative equation was determined, i.e., the optimized values of the input parameters $\alpha$, $S(x,t)$, $k$, $h$, $T_\infty$. For the coarse optimization the first method used was Conjugate Gradient (CG). The accuracy(how close the answer of fine model and coarse model are) of 89.2\% was achieved. The number of hidden layers for the neural network was 3. The project was assessed on another optimization method for the coarse model and that is Nelder-Mead method. The accuracy achieved with that was 90.1\%. 

Another observation was that, when the number of hidden layers(for the neural network) increases the accuracy was increasing upto some level. When the number of hidden layers was less than 3 , the neural network went through underfitting and the accuracy fell to 60.2\%. On the otherhand, when the hidden layers were chosen to be more than 6, the model became too complex and went through overfitting and that made the accuracy to fall to approximately 65\%. Hence a proper consideration has to be given to the choice of the number of hidden layers and the optimization method to be used for the coarse model.

\subsection{Future work}
\vspace{3pt}
The current implementation uses the trapezoidal rule to integrate the differential equation in time. However, other numerical methods like the Runge-Kutta\cite{abc} method or finite difference methods[like forward or backward process] to obtain more accurate solutions. To evaluate the generic behaviour of the proposed algorithm, it can be checked on more complex PDEs like Schrodinger equation, which models quantum mechanics. 

Also, the training data in the current implementation was generated by adding Gaussian noise to the analytical solution. So, the effect of different levels of noise on the accuracy of the neural network model can be investigated as a future work. Morover, the proposed method can be checked on other optimization techniques applied on the coarse model. The other thing, which seems exciting, and could be done, is the use of other machine learning models like support vector machines (SVMs), random forests, or convolutional neural networks (CNNs) to solve the PDE, it sounds challenging but really worth doing.

\printbibliography

\end{document}